# CONVERGENCE RATE OF LINEAR TWO-TIME-SCALE STOCHASTIC APPROXIMATION[1]

By Vijay R. Konda and John N. Tsitsiklis

*Massachusetts Institute of Technology*

We study the rate of convergence of linear two-time-scale stochastic approximation methods. We consider two-time-scale linear iterations driven by i.i.d. noise, prove some results on their asymptotic covariance and establish asymptotic normality. The well-known result [Polyak, B. T. (1990). *Automat. Remote Contr.* **51** 937–946; Ruppert, D. (1988). Technical Report 781, Cornell Univ.] on the optimality of Polyak–Ruppert averaging techniques specialized to linear stochastic approximation is established as a consequence of the general results in this paper.

**1. Introduction.** Two-time-scale stochastic approximation methods [Borkar (1997)] are recursive algorithms in which some of the components are updated using step-sizes that are very small compared to those of the remaining components. Over the past few years, several such algorithms have been proposed for various applications [Konda and Borkar (1999), Bhatnagar, Fu, Marcus and Fard (2001), Baras and Borkar (2000), Bhatnagar, Fu and Marcus (2001) and Konda and Tsitsiklis (2003)].

The general setting for two-time-scale algorithms is as follows. Let $f(\theta, r)$ and $g(\theta, r)$ be two unknown functions and let $(\theta^*, r^*)$ be the unique solution to the equations

$$(1.1) \qquad f(\theta, r) = 0, \qquad g(\theta, r) = 0.$$

The functions $f(\cdot, \cdot)$ and $g(\cdot, \cdot)$ are accessible only by simulating or observing a stochastic system which, given $\theta$ and $r$ as input, produces $F(\theta, r, V)$ and $G(\theta, r, W)$. Here, $V$ and $W$ are random variables, representing noise, whose distribution satisfies

$$f(\theta, r) = E[F(\theta, r, V)], \qquad g(\theta, r) = E[G(\theta, r, W)] \qquad \forall\, \theta, r.$$

Received March 2002; revised March 2003.
[1]Supported by NSF Grant ECS-9873451.
*AMS 2000 subject classification.* 62L20.
*Key words and phrases.* Stochastic approximation, two-time-scales.

Assume that the noise $(V, W)$ in each simulation or observation of the stochastic system is independent of the noise in all other simulations. In other words, assume that we have access to an independent sequence of functions $F(\cdot, \cdot, V_k)$ and $G(\cdot, \cdot, W_k)$. Suppose that for any given $\theta$, the stochastic iteration

$$r_{k+1} = r_k + \gamma_k G(\theta, r_k, W_k) \tag{1.2}$$

is known to converge to some $h(\theta)$. Furthermore, assume that the stochastic iteration

$$\theta_{k+1} = \theta_k + \gamma_k F(\theta_k, h(\theta_k), V_k) \tag{1.3}$$

is known to converge to $\theta^*$. Given this information, we wish to construct an algorithm that solves the system of equations (1.1).

Note that the iteration (1.2) has only been assumed to converge when $\theta$ is held fixed. This assumption allows us to fix $\theta$ at a current value $\theta_k$, run the iteration (1.2) for a long time, so that $r_k$ becomes approximately equal to $h(\theta_k)$, use the resulting $r_k$ to update $\theta_k$ in the direction of $F(\theta_k, r_k, W_k)$, and repeat this procedure. While this is a sound approach, it requires an increasingly large time between successive updates of $\theta_k$. Two-time-scale stochastic approximation methods circumvent this difficulty by using different step sizes $\{\beta_k\}$ and $\{\gamma_k\}$ and update $\theta_k$ and $r_k$, according to

$$\theta_{k+1} = \theta_k + \beta_k F(\theta_k, r_k, V_k),$$
$$r_{k+1} = r_k + \gamma_k G(\theta_k, r_k, W_k),$$

where $\beta_k$ is very small relative to $\gamma_k$. This makes $\theta_k$ "quasi-static" compared to $r_k$ and has an effect similar to fixing $\theta_k$ and running the iteration (1.2) forever. In turn, $\theta_k$ sees $r_k$ as a close approximation of $h(\theta_k)$ and therefore its update looks almost the same as (1.3).

How small should the ratio $\beta_k/\gamma_k$ be for the above scheme to work? The answer generally depends on the functions $f(\cdot, \cdot)$ and $g(\cdot, \cdot)$, which are typically unknown. This leads us to consider a safe choice whereby $\beta_k/\gamma_k \to 0$. The subject of this paper is the convergence rate analysis of the two-time-scale algorithms that result from this choice. We note here that the analysis is significantly different from the case where $\lim_k(\beta_k/\gamma_k) > 0$, which can be handled using existing techniques.

Two-time-scale algorithms have been proved to converge in a variety of contexts [Borkar (1997), Konda and Borkar (1999) and Konda and Tsitsiklis (2003)]. However, except for the special case of Polyak–Ruppert averaging, there are no results on their rate of convergence. The existing analysis [Ruppert (1988), Polyak (1990), Polyak and Juditsky (1992) and Kushner and Yang (993)] of Polyak–Ruppert methods rely on special structure and are not applicable to the more general two-time-scale iterations considered here.

The main result of this paper is a rule of thumb for calculating the asymptotic covariance of linear two-time-scale stochastic iterations. For example, consider the linear iterations

$$\theta_{k+1} = \theta_k + \beta_k(b_1 - A_{11}\theta_k - A_{12}r_k + V_k), \tag{1.4}$$

$$r_{k+1} = r_k + \gamma_k(b_2 - A_{21}\theta_k - A_{22}r_k + W_k). \tag{1.5}$$

We show that the asymptotic covariance matrix of $\beta_k^{-1/2}\theta_k$ is the same as that of $\beta_k^{-1/2}\bar{\theta}_k$, where $\bar{\theta}_k$ evolves according to the single-time-scale stochastic iteration:

$$\bar{\theta}_{k+1} = \bar{\theta}_k + \beta_k(b_1 - A_{11}\bar{\theta}_k - A_{12}\bar{r}_k + V_k),$$
$$0 = b_2 - A_{21}\bar{\theta}_k - A_{22}\bar{r}_k + W_k.$$

Besides the calculation of the asymptotic covariance of $\beta_k^{-1/2}\theta_k$ (Theorem 2.8), we also establish that the distribution of $\beta_k^{-1/2}(\theta_k - \theta^*)$ converges to a Gaussian with mean zero and with the above asymptotic covariance (Theorem 4.1). We believe that the proof techniques of this paper can be extended to nonlinear stochastic approximation to obtain similar results. However, this and other possible extensions (such as weak convergence of paths to a diffusion process) are no pursued in this paper.

In the linear case, our results also explain why Polyak–Ruppert averaging is optimal. Suppose that we are looking for the solution of the linear system

$$Ar = b$$

in a setting where we only have access to noisy measurements of $b - Ar$. The standard algorithm in this setting is

$$r_{k+1} = r_k + \gamma_k(b - Ar_k + W_k), \tag{1.6}$$

and is known to converge under suitable conditions. (Here, $W_k$ represents zero-mean noise at time $k$.) In order to improve the rate of convergence, Polyak (1990) and Ruppert (1988) suggest using the average

$$\theta_k = \frac{1}{k}\sum_{l=0}^{k-1} r_l \tag{1.7}$$

as an estimate of the solution, instead of $r_k$. It was shown in Polyak (1990) that if $k\gamma_k \to \infty$, the asymptotic covariance of $\sqrt{k}\theta_k$ is $A^{-1}\Gamma(A')^{-1}$, where $\Gamma$ is the covariance of $W_k$. Furthermore, this asymptotic covariance matrix is known to be optimal [Kushner and Yin (1997)].

The calculation of the asymptotic covariance in Polyak (1990) and Ruppert (1988) uses the special averaging structure. We provide here an alternative calculation based on our results. Note that $\theta_k$ satisfies the recursion

$$\theta_{k+1} = \theta_k + \frac{1}{k+1}(r_k - \theta_k), \tag{1.8}$$

and the iteration (1.6)–(1.8) for $r_k$ and $\theta_k$ is a special case of the two-time-scale iterations (1.4) and (1.5), with the correspondence $b_1 = 0$, $A_{11} = I$, $A_{12} = -I$, $V_k = 0$, $b_2 = b$, $A_{21} = 0$, $A_{22} = 0$. Furthermore, the assumption $k\gamma_k \to \infty$ corresponds to our general assumption $\beta_k/\gamma_k \to 0$.

By applying our rule of thumb to the iteration (1.6)–(1.8), we see that the asymptotic covariance of $(\sqrt{k+1})\theta_k$ is the same as that of $(\sqrt{k+1})\bar{\theta}_k$, where $\bar{\theta}_k$ satisfies

$$\bar{\theta}_{k+1} = \bar{\theta}_k + \frac{1}{k+1}(-\bar{\theta}_k + A^{-1}(b + W_k)),$$

or

$$\bar{\theta}_k = \frac{1}{k}\sum_{l=0}^{k-1}(A^{-1}b + A^{-1}W_l).$$

It then follows that the covariance of $\sqrt{k}\bar{\theta}_k$ is $A^{-1}\Gamma(A')^{-1}$, and we recover the result of Polyak (1990), Polyak and Juditsky (1992) and Ruppert (1988) for the linear case.

In the example just discussed, the use of two time-scales is not necessary for convergence, but is essential for the improvement of the convergence rate. This idea of introducing two time-scales to improve the rate of convergence deserves further exploration. It is investigated to some extent in the context of reinforcement learning algorithms in Konda (2002).

Finally, we would like to point out the differences between the two-time-scale iterations we study here and those that arise in the study of the tracking ability of adaptive algorithms [see Benveniste, Metivier and Priouret (1990)]. There, the slow component represents the movement of underlying system parameters and the fast component represents the user's algorithm. The fast component, that is, the user's algorithm, does not affect the slow component. In contrast, we consider iterations in which the fast component affects the slow one and vice versa. Furthermore, the relevant figures of merit are different. For example, in Benveniste, Metivier and Priouret (1990), one is mostly interested in the behavior of the fast component, whereas we focus on the asymptotic covariance of the slow component.

The outline of the paper is as follows. In the next section, we consider linear iterations driven by i.i.d. noise and obtain expressions for the asymptotic covariance of the iterates. In Section 3, we compare the convergence rate of two-time-scale algorithms and their single-time-scale counterparts. In Section 4, we establish asymptotic normality of the iterates.

Before proceeding, we introduce some notation. Throughout the paper, $|\cdot|$ represents the Euclidean norm of vectors or the induced operator norm of matrices. Furthermore, $I$ and $0$ represent identity and null matrices, respectively. We use the abbreviation w.p.1 for "with probability 1." We use $c, c_1, c_2, \ldots$ to represent some constants whose values are not important.

**2. Linear iterations.** In this section, we consider iterations of the form

(2.1) $$\theta_{k+1} = \theta_k + \beta_k(b_1 - A_{11}\theta_k - A_{12}r_k + V_k),$$

(2.2) $$r_{k+1} = r_k + \gamma_k(b_2 - A_{21}\theta_k - A_{22}r_k + W_k),$$

where $\theta_k$ is in $\mathbf{R}^n$, $r_k$ is in $\mathbf{R}^m$, and $b_1$, $b_2$, $A_{11}$, $A_{12}$, $A_{21}$, $A_{22}$ are vectors and matrices of appropriate dimensions.

Before we present our results, we motivate various assumptions that we will need. The first two assumptions are standard.

ASSUMPTION 2.1. The random variables $(V_k, W_k)$, $k = 0, 1, \ldots$, are independent of $r_0$, $\theta_0$, and of each other. They have zero mean and common covariance

$$E[V_k V_k'] = \Gamma_{11},$$
$$E[V_k W_k'] = \Gamma_{12} = \Gamma_{21}',$$
$$E[W_k W_k'] = \Gamma_{22}.$$

ASSUMPTION 2.2. The step-size sequences $\{\gamma_k\}$ and $\{\beta_k\}$ are deterministic, positive, nonincreasing, and satisfy the following:

1. $\sum_k \gamma_k = \sum_k \beta_k = \infty$.
2. $\beta_k, \gamma_k \to 0$.

The key assumption that the step sizes $\beta_k$ and $\gamma_k$ are of different orders of magnitude is subsumed by the following.

ASSUMPTION 2.3. There exists some $\varepsilon \geq 0$ such that

$$\frac{\beta_k}{\gamma_k} \to \varepsilon.$$

For the iterations (2.1) and (2.2) to be consistent with the general scheme of two-time-scale stochastic approximations described in the Introduction, we need some assumptions on the matrices $A_{ij}$. In particular, we need iteration (2.2) to converge to $A_{22}^{-1}(b_2 - A_{21}\theta)$, when $\theta_k$ is held constant at $\theta$. Furthermore, the sequence $\theta_k$ generated by the iteration

$$\theta_{k+1} = \theta_k + \beta_k(b_1 - A_{12}A_{22}^{-1}b_2 - (A_{11} - A_{12}A_{22}^{-1}A_{21})\theta_k + V_k),$$

which is obtained by substituting $A_{22}^{-1}(b_2 - A_{21}\theta_k)$ for $r_k$ in iteration (2.1), should also converge. Our next assumption is needed for the above convergence to take place.

Let $\Delta$ be the matrix defined by

(2.3) $$\Delta = A_{11} - A_{12}A_{22}^{-1}A_{21}.$$

Recall that a square matrix $A$ is said to be Hurwitz if the real part of each eigenvalue of $A$ is strictly negative.

ASSUMPTION 2.4. The matrices $-A_{22}, -\Delta$ are Hurwitz.

It is not difficult to show that, under the above assumptions, $(\theta_k, r_k)$ converges in mean square and w.p.1 to $(\theta^*, r^*)$. The objective of this paper is to capture the rate at which this convergence takes place. Obviously, this rate depends on the step-sizes $\beta_k, \gamma_k$, and this dependence can be quite complicated in general. The following assumption ensures that the rate of mean square convergence of $(\theta_k, r_k)$ to $(\theta^*, r^*)$ bears a simple relationship (asymptotically linear) with the step-sizes $\beta_k, \gamma_k$.

ASSUMPTION 2.5.   1. There exists a constant $\bar{\beta} \geq 0$ such that
$$\lim_k (\beta_{k+1}^{-1} - \beta_k^{-1}) = \bar{\beta}.$$

2. If $\varepsilon = 0$, then
$$\lim_k (\gamma_{k+1}^{-1} - \gamma_k^{-1}) = 0.$$

3. The matrix $-(\Delta - \frac{\bar{\beta}}{2}I)$ is Hurwitz.

Note that when $\varepsilon > 0$, the iterations (2.1) and (2.2) are essentially single-time-scale algorithms and therefore can be analyzed using existing techniques [Nevel'son and Has'minskii (1973), Kusher and Clark (1978), Benveniste,
Metivier and Priouret (1990), Duflo (1997) and Kusher and Yin (1997)]. We include this in our analysis as we would like to study the behavior of the rate of convergence as $\varepsilon \downarrow 0$. The following is an example of sequences satisfying the above assumption with $\varepsilon = 0$, $\bar{\beta} = 1/(\tau_1 \beta_0)$:

$$\gamma_k = \frac{\gamma_0}{(1 + k/\tau_0)^\alpha}, \qquad \frac{1}{2} < \alpha < 1,$$

$$\beta_k = \frac{\beta_0}{(1 + k/\tau_1)},$$

Let $\theta^* \in \mathbf{R}^m$ and $r^* \in \mathbf{R}^n$ be the unique solution to the system of linear equations

$$A_{11}\theta + A_{12}r = b_1,$$
$$A_{21}\theta + A_{22}r = b_2.$$

For each $k$, let

$$\hat{\theta}_k = \theta_k - \theta^*,$$
(2.4)
$$\hat{r}_k = r_k - A_{22}^{-1}(b_2 - A_{21}\theta_k)$$

and

$$\Sigma_{11}^k = \beta_k^{-1} E[\hat{\theta}_k \hat{\theta}_k'],$$
$$\Sigma_{12}^k = (\Sigma_{21}^k)' = \beta_k^{-1} E[\hat{\theta}_k \hat{r}_k'],$$
$$\Sigma_{22}^k = \gamma_k^{-1} E[\hat{r}_k \hat{r}_k'],$$
$$\Sigma^k = \begin{bmatrix} \Sigma_{11}^k & \Sigma_{12}^k \\ \Sigma_{21}^k & \Sigma_{22}^k \end{bmatrix}.$$

Our main result is the following.

THEOREM 2.6. *Under Assumptions 2.1–2.5, and when the constant $\varepsilon$ of Assumption 2.3 is sufficiently small, the limit matrices*

(2.5)  $\Sigma_{11}^{(\varepsilon)} = \lim_k \Sigma_{11}^k, \qquad \Sigma_{12}^{(\varepsilon)} = \lim_k \Sigma_{12}^k, \qquad \Sigma_{22}^{(\varepsilon)} = \lim_k \Sigma_{22}^k$

*exist. Furthermore, the matrix*

$$\Sigma^{(0)} = \begin{bmatrix} \Sigma_{11}^{(0)} & \Sigma_{12}^{(0)} \\ \Sigma_{21}^{(0)} & \Sigma_{22}^{(0)} \end{bmatrix}$$

*is the unique solution to the following system of equations*

(2.6) $\quad \Delta \Sigma_{11}^{(0)} + \Sigma_{11}^{(0)} \Delta' - \bar{\beta} \Sigma_{11}^{(0)} + A_{12} \Sigma_{21}^{(0)} + \Sigma_{12}^{(0)} A_{12}' = \Gamma_{11},$

(2.7) $\hspace{10em} A_{12} \Sigma_{22}^{(0)} + \Sigma_{12}^{(0)} A_{22}' = \Gamma_{12},$

(2.8) $\hspace{10em} A_{22} \Sigma_{22}^{(0)} + \Sigma_{22}^{(0)} A_{22}' = \Gamma_{22}.$

*Finally,*

(2.9) $\quad \lim_{\varepsilon \downarrow 0} \Sigma_{11}^{(\varepsilon)} = \Sigma_{11}^{(0)}, \qquad \lim_{\varepsilon \downarrow 0} \Sigma_{12}^{(\varepsilon)} = \Sigma_{12}^{(0)}, \qquad \lim_{\varepsilon \downarrow 0} \Sigma_{22}^{(\varepsilon)} = \Sigma_{22}^{(0)}.$

PROOF. Let us first consider the case $\varepsilon = 0$. The idea of the proof is to study the iteration in terms of transformed variables:

(2.10) $\hspace{8em} \tilde{\theta}_k = \hat{\theta}_k, \qquad \tilde{r}_k = L_k \hat{\theta}_k + \hat{r}_k,$

for some sequence of $n \times m$ matrices $\{L_k\}$ which we will choose so that *the faster time-scale iteration does not involve the slower time-scale variables.* To see what the sequence $\{L_k\}$ should be, we rewrite the iterations (2.1) and (2.2) in terms of the transformed variables as shown below (see Section A.1 for the algebra leading to these equations):

$$
\begin{aligned}
\tilde{\theta}_{k+1} &= \tilde{\theta}_k - \beta_k(B_{11}^k \tilde{\theta}_k + A_{12}\tilde{r}_k) + \beta_k V_k, \\
\tilde{r}_{k+1} &= \tilde{r}_k - \gamma_k(B_{21}^k \tilde{\theta}_k + B_{22}^k \tilde{r}_k) + \gamma_k W_k + \beta_k(L_{k+1} + A_{22}^{-1}A_{21})V_k,
\end{aligned}
\tag{2.11}
$$

where

$$
\begin{aligned}
B_{11}^k &= \Delta - A_{12}L_k, \\
B_{21}^k &= \frac{L_k - L_{k+1}}{\gamma_k} + \frac{\beta_k}{\gamma_k}(L_{k+1} + A_{22}^{-1}A_{21})B_{11}^k - A_{22}L_k, \\
B_{22}^k &= \frac{\beta_k}{\gamma_k}(L_{k+1} + A_{22}^{-1}A_{21})A_{12} + A_{22}.
\end{aligned}
$$

We wish to choose $\{L_k\}$ so that $B_{21}^k$ is eventually zero. To accomplish this, we define the sequence of matrices $\{L_k\}$ by

$$
\begin{aligned}
L_k &= 0, \qquad 0 \leq k \leq k_0, \\
L_{k+1} &= (L_k - \gamma_k A_{22}L_k + \beta_k A_{22}^{-1}A_{21}B_{11}^k)(I - \beta_k B_{11}^k)^{-1} \qquad \forall k \geq k_0,
\end{aligned}
\tag{2.12}
$$

so that $B_{21}^k = 0$ for all $k \geq k_0$. For the above recursion to be meaningful, we need $(I - \beta_k B_{11}^k)$ to be nonsingular for all $k \geq k_0$. This is handled by Lemma A.1 in the Appendix, which shows that if $k_0$ is sufficiently large, then the sequence of matrices $\{L_k\}$ is well defined and also converges to zero.

For every $k \geq k_0$, we define

$$
\begin{aligned}
\tilde{\Sigma}_{11}^k &= \beta_k^{-1} E[\tilde{\theta}_k \tilde{\theta}_k'], \\
(\tilde{\Sigma}_{21}^k)' &= \tilde{\Sigma}_{12}^k = \beta_k^{-1} E[\tilde{\theta}_k \tilde{r}_k'], \\
\tilde{\Sigma}_{22}^k &= \gamma_k^{-1} E[\tilde{r}_k \tilde{r}_k'].
\end{aligned}
$$

Using the transformation (2.10), it is easy to see that

$$
\begin{aligned}
\tilde{\Sigma}_{11}^k &= \Sigma_{11}^k, \\
\tilde{\Sigma}_{12}^k &= \Sigma_{11}^k L_k' + \Sigma_{12}^k, \\
\tilde{\Sigma}_{22}^k &= \Sigma_{22}^k + \left(\frac{\beta_k}{\gamma_k}\right)(L_k \Sigma_{12}^k + \Sigma_{21}^k L_k' + L_k \Sigma_{11}^k L_k').
\end{aligned}
$$

Since $L_k \to 0$, we obtain

$$
\begin{aligned}
\lim_k \Sigma_{11}^k &= \lim_k \tilde{\Sigma}_{11}^k, \\
\lim_k \Sigma_{12}^k &= \lim_k \tilde{\Sigma}_{12}^k, \\
\lim_k \Sigma_{22}^k &= \lim_k \tilde{\Sigma}_{12}^k,
\end{aligned}
$$

provided that the limits exist.

To compute $\lim_k \tilde{\Sigma}_{22}^k$, we use (2.11), the fact that $B_{21}^k = 0$ for large enough $k$, the fact that $B_{22}^k$ converges to $A_{22}$, and some algebra, to arrive at the following recursion for $\tilde{\Sigma}_{22}^k$:

$$(2.13) \quad \tilde{\Sigma}_{22}^{k+1} = \tilde{\Sigma}_{22}^k + \gamma_k(\Gamma_{22} - A_{22}\tilde{\Sigma}_{22}^k - \tilde{\Sigma}_{22}^k A_{22}' + \delta_{22}^k(\tilde{\Sigma}_{22}^k)),$$

where $\delta_{22}^k(\cdot)$ is some matrix-valued affine function (on the space of matrices) such that

$$\lim_k \delta_{22}^k(\Sigma_{22}) = 0 \qquad \text{for all } \Sigma_{22}.$$

Since $-A_{22}$ is Hurwitz, it follows (see Lemma A.2 in the Appendix) that the limit

$$\lim_k \Sigma_{22}^k = \lim_k \tilde{\Sigma}_{22}^k = \Sigma_{22}^{(0)}$$

exists, and $\Sigma_{22}^{(0)}$ satisfies (2.8).

Similarly, $\tilde{\Sigma}_{12}^k$ satisfies

$$(2.14) \quad \tilde{\Sigma}_{12}^{k+1} = \tilde{\Sigma}_{12}^k + \gamma_k(\Gamma_{12} - A_{12}\Sigma_{22}^{(0)} - \tilde{\Sigma}_{12}^k A_{22}' + \delta_{12}^k(\tilde{\Sigma}_{12}^k))$$

where, as before, $\delta_{12}^k(\cdot)$ is an affine function that goes to zero. (The coefficients of this affine function depend, in general, on $\tilde{\Sigma}_{22}^k$, but the important property is that they tend to zero as $k \to \infty$.) Since $-A_{22}$ is Hurwitz, the limit

$$\lim_k \Sigma_{12}^k = \lim_k \tilde{\Sigma}_{12}^k = \Sigma_{12}^{(0)}$$

exists and satisfies (2.7). Finally, $\tilde{\Sigma}_{11}^k$ satisfies

$$(2.15) \quad \begin{aligned}\tilde{\Sigma}_{11}^{k+1} = \tilde{\Sigma}_{11}^k + \beta_k(\Gamma_{11} - A_{12}\Sigma_{21}^{(0)} - \Sigma_{12}^{(0)} A_{12}' - \Delta\tilde{\Sigma}_{11}^k \\ - \tilde{\Sigma}_{11}^k \Delta' + \bar{\beta}\tilde{\Sigma}_{11}^k + \delta_{11}^k(\tilde{\Sigma}_{11}^k)),\end{aligned}$$

where $\delta_{11}^k(\cdot)$ is some affine function that goes to zero. (Once more, the coefficients of this affine function depend, in general, on $\tilde{\Sigma}_{22}^k$ and $\tilde{\Sigma}_{12}^k$, but they tend to zero as $k \to \infty$.) Since $-(\Delta - \frac{\bar{\beta}}{2}I)$ is Hurwitz, the limit

$$\lim_k \Sigma_{11}^k = \lim_k \tilde{\Sigma}_{11}^k = \Sigma_{11}^{(0)}$$

exists and satisfies (2.6).

The above arguments show that for $\varepsilon = 0$, the limit matrices in (2.5) exist and satisfy (2.6)–(2.8). To complete the proof, we need to show that these limit matrices exist for sufficiently small $\varepsilon > 0$ and that the limiting relations (2.9) hold. As this part of the proof uses standard techniques, we will only outline the analysis.

Define for each $k$,

$$Z_k = \begin{pmatrix} \hat{\theta}_k \\ \hat{r}_k \end{pmatrix}.$$

The linear iterations (2.1) and (2.2) can be rewritten in terms of $Z_k$ as

$$Z_{k+1} = Z_k - \beta_k B_k Z_k + \beta_k U_k,$$

where $U_k$ is a sequence of independent random vectors and $\{B_k\}$ is a sequence of deterministic matrices. Using the assumption that $\beta_k/\gamma_k$ converges to $\varepsilon$, it can be shown that the sequence of matrices $B_k$ converges to some matrix $B^{(\varepsilon)}$ and, similarly, that

$$\lim_k E[U_k U_k'] = \Gamma^{(\varepsilon)}$$

for some matrix $\Gamma^{(\varepsilon)}$. Furthermore, when $\varepsilon > 0$ is sufficiently small, it can be shown that $-(B^{(\varepsilon)} - \frac{\bar{\beta}}{2}I)$ is Hurwitz. It then follows from standard theorems [see, e.g., Polyak (1976)] on the asymptotic covariance of stochastic approximation methods, that the limit

$$\lim_k \beta_k^{-1} E[Z_k Z_k']$$

exists and satisfies a *linear* equation whose coefficients depend smoothly on $\varepsilon$ (the coefficients are infinitely differentiable w.r.t. $\varepsilon$). Since the components of the above limit matrix are $\Sigma_{11}^{(\varepsilon)}$, $\Sigma_{12}^{(\varepsilon)}$ and $\Sigma_{22}^{(\varepsilon)}$ modulo some scaling, the latter matrices also satisfy a linear equation which depends on $\varepsilon$. The explicit form of this equation is tedious to write down and does not provide any additional insight for our purposes. We note, however, that when we set $\varepsilon$ to zero, this system of equations becomes the same as (2.6)–(2.8). Since (2.6)–(2.8) have a unique solution, the system of equations for $\Sigma_{11}^{(\varepsilon)}$, $\Sigma_{12}^{(\varepsilon)}$ and $\Sigma_{22}^{(\varepsilon)}$ also has a unique solution for all sufficiently small $\varepsilon$. Furthermore, the dependence of the solution on $\varepsilon$ is smooth because the coefficients are smooth in $\varepsilon$. □

REMARK 2.7. The transformations used in the above proof are inspired by those used to study singularly perturbed ordinary differential equations [Kokotovic (1984)]. However, most of these transformations were time-invariant because the perturbation parameter was constant. In such cases, the matrix $L$ satisfies a static Riccati equation instead of the recursion (2.12). In contrast, our transformations are time-varying because our "perturbation" parameter $\beta_k/\gamma_k$ is time-varying.

In most applications, the iterate $r_k$ corresponds to some auxiliary parameters and one is mostly interested in the asymptotic covariance $\Sigma_{11}^{(0)}$ of $\theta_k$.

Note that according to Theorem 2.6, the covariance of the auxiliary parameters is of the order of $\gamma_k$, whereas the covariance of $\theta_k$ is of the order of $\beta_k$. With two time-scales, one can potentially improve the rate of convergence of $\theta_k$ (cf. to a single-time-scale algorithm) by sacrificing the rate of convergence of the auxiliary parameters. To make such comparisons possible, we need an alternative interpretation of $\Sigma_{11}^{(0)}$, that does not explicitly refer to the system (2.6)–(2.8). This is accomplished by our next result, which provides a useful tool for the design and analysis of two-time-scale stochastic approximation methods.

THEOREM 2.8. *The asymptotic covariance matrix $\Sigma_{11}^{(0)}$ of $\beta_k^{-1/2}\theta_k$ is the same as the asymptotic covariance of $\beta_k^{-1/2}\bar{\theta}_k$, where $\bar{\theta}_k$ is generated by*

$$(2.16) \qquad \bar{\theta}_{k+1} = \bar{\theta}_k + \beta_k(b_1 - A_{11}\bar{\theta}_k - A_{12}\bar{r}_k + V_k),$$

$$(2.17) \qquad 0 = b_2 - A_{21}\bar{\theta}_k - A_{22}\bar{r}_k + W_k.$$

*In other words,*

$$\Sigma_{11}^{(0)} = \lim_k \beta_k^{-1} E[\bar{\theta}_k \bar{\theta}_k'].$$

PROOF. We start with (2.6)–(2.8) and perform some algebraic manipulations to eliminate $\Sigma_{12}^{(0)}$ and $\Sigma_{22}^{(0)}$. This leads to a single equation for $\Sigma_{11}^{(0)}$, of the form

$$\Delta \Sigma_{11}^{(0)} + \Sigma_{11}^{(0)} \Delta' - \bar{\beta} \Sigma_{11}^{(0)}$$
$$= \Gamma_{11} - A_{12} A_{22}^{-1} \Gamma_{21} - \Gamma_{12} (A_{22}')^{-1} A_{12}' + A_{12} A_{22}^{-1} \Gamma_{22} (A_{22}')^{-1} A_{12}'.$$

Note that the right-hand side of the above equation is exactly the covariance of $V_k - A_{12} A_{22}^{-1} W_k$. Therefore, the asymptotic covariance of $\theta_k$ is the same as the asymptotic covariance of the following stochastic approximation:

$$\bar{\theta}_{k+1} = \bar{\theta}_k + \beta_k(-\Delta \bar{\theta}_k + V_k - A_{12} A_{22}^{-1} W_k).$$

Finally, note that the above iteration is the one obtained by eliminating $r_k$ from iterations (2.16) and (2.17). □

REMARK. The single-time-scale stochastic approximation procedure in Theorem 2.8 is not implementable when the matrices $A_{ij}$ are unknown. The theorem establishes that two-time-scale stochastic approximation performs as well as if these matrices are known.

REMARK. The results of the previous section show that the asymptotic covariance matrix of $\beta_k^{-1/2}\theta_k$ is independent of the step-size schedule $\{\gamma_k\}$ for the fast iteration if

$$\frac{\beta_k}{\gamma_k} \to 0.$$

To understand, at least qualitatively, the effect of the step-sizes $\gamma_k$ on the transient behavior, recall the recursions (2.13)–(2.15) satisfied by the covariance matrices $\tilde{\Sigma}^k$:

$$\tilde{\Sigma}_{11}^{k+1} = \tilde{\Sigma}_{11}^k + \beta_k(\Gamma_{11} - A_{12}\Sigma_{21}^{(0)} - \Sigma_{12}^{(0)}A_{12}'$$
$$- \Delta\tilde{\Sigma}_{11}^k - \tilde{\Sigma}_{11}^k\Delta' - \bar{\beta}\tilde{\Sigma}_{11}^k + \delta_{11}^k(\tilde{\Sigma}_{11}^k)),$$
$$\tilde{\Sigma}_{12}^{k+1} = \tilde{\Sigma}_{12}^k + \gamma_k(\Gamma_{12} - A_{12}\Sigma_{22}^{(0)} - \tilde{\Sigma}_{12}^k A_{22}' + \delta_{12}^k(\tilde{\Sigma}_{12}^k)),$$
$$\tilde{\Sigma}_{22}^{k+1} = \Sigma_{22}^k + \gamma_k(\Gamma_{22} - A_{22}\Sigma_{22}^k - \Sigma_{22}^k A_{22}' + \delta_{22}^k(\Sigma_{22}^k)),$$

where the $\delta_{ij}^k(\cdot)$ are affine functions that tend to zero as $k$ tends to infinity. Using explicit calculations, it is easy to verify that the error terms $\delta_{ij}^k$ are of the form

$$\delta_{11}^k = A_{12}(\tilde{\Sigma}_{21}^k - \Sigma_{21}^{(0)}) + (\tilde{\Sigma}_{12}^k - \Sigma_{12}^{(0)})A_{12}' + O(\beta_k),$$
$$\delta_{12}^k = A_{12}(\Sigma_{22}^{(0)} - \tilde{\Sigma}_{22}^k) + O\left(\frac{\beta_k}{\gamma_k}\right),$$
$$\delta_{22}^k = O\left(\frac{\beta_k}{\gamma_k}\right).$$

To clarify the meaning of the above relations, the first one states that the affine function $\delta_{11}^k(\Sigma_{11})$ is the sum of the constant term $A_{12}(\tilde{\Sigma}_{21}^k - \Sigma_{21}^{(0)}) + (\tilde{\Sigma}_{12}^k - \Sigma_{12}^{(0)})A_{12}'$, and another affine function of $\Sigma_{11}^k$ whose coefficients are proportional to $\beta_k$.

The above relations show that the rate at which $\tilde{\Sigma}_{11}^k$ converges to $\Sigma_{11}^{(0)}$ depends on the rate at which $\tilde{\Sigma}_{12}^k$ converges to $\Sigma_{12}^{(0)}$, through the term $\delta_{11}^k$. The rate of convergence of $\tilde{\Sigma}_{12}^k$, in turn, depends on that of $\tilde{\Sigma}_{22}^k$, through the term $\delta_{12}^k$. Since the step-size in the recursions for $\tilde{\Sigma}_{22}^k$ and $\tilde{\Sigma}_{12}^k$ is $\gamma_k$, and the error terms in these recursions are proportional to $\beta_k/\gamma_k$, the transients depend on both sequences $\{\gamma_k\}$ and $\{\beta_k/\gamma_k\}$. But each sequence has a different effect. When $\gamma_k$ is large, instability or large oscillations of $r_k$ are possible. On the other hand, when $\beta_k/\gamma_k$ is large, the error terms $\delta_{ij}^k$ can be large and can prolong the transient period. Therefore, one would like to have $\beta_k/\gamma_k$ decrease to zero quickly, while at the same time avoiding large $\gamma_k$. Apart from these loose guidelines, it appears difficult to obtain a characterization of desirable step-size schedules.

**3. Single time-scale versus two time-scales.** In this section, we compare the optimal asymptotic covariance of $\beta_k^{-1/2}\theta_k$ that can be obtained by a realizable single-time-scale stochastic iteration, with the optimal asymptotic covariance that can be obtained by a realizable two-time-scale stochastic iteration. The optimization is to be carried out over a set of suitable gain matrices that can be used to modify the algorithm, and the optimality criterion to be used is one whereby a covariance matrix $\Sigma$ is preferable to another covariance matrix $\tilde{\Sigma}$ if $\tilde{\Sigma} - \Sigma$ is nonzero and nonnegative definite.

Recall that Theorem 2.8 established that the asymptotic covariance of a two-time-scale iteration is the same as in a related single-time-scale iteration. However, the related single-time-scale iteration is unrealizable, unless the matrix $A$ is known. In contrast, in this section we compare realizable iterations that do not require explicit knowledge of $A$ (although knowledge of $A$ would be required in order to select the best possible realizable iteration).

We now specify the classes of stochastic iterations that we will be comparing.

1. We consider two-time-scale iterations of the form

$$\theta_{k+1} = \theta_k + \beta_k G_1 (b_1 - A_{11}\theta_k - A_{12}r_k + V_k),$$
$$r_{k+1} = r_k + \gamma_k (b_2 - A_{21}\theta_k - A_{22}r_k + W_k).$$

Here, $G_1$ is a gain matrix, which we are allowed to choose in a manner that minimizes the asymptotic covariance of $\beta_k^{-1/2}\theta_k$.

2. We consider single-time-scale iterations, in which we have $\gamma_k = \beta_k$, but in which we are allowed to use an arbitrary gain matrix $G$, in order to minimize the asymptotic covariance of $\beta_k^{-1/2}\theta_k$. Concretely, we consider iterations of the form

$$\begin{bmatrix} \theta_{k+1} \\ r_{k+1} \end{bmatrix} = \begin{bmatrix} \theta_k \\ r_k \end{bmatrix} + \beta_k G \begin{bmatrix} b_1 - A_{11}\theta_k - A_{12}r_k + V_k \\ b_2 - A_{21}\theta_k - A_{22}r_k + W_k \end{bmatrix}.$$

We then have the following result.

THEOREM 3.1. *Under Assumptions 2.1–2.5, and with $\varepsilon = 0$, the minimal possible asymptotic covariance of $\beta_k^{-1/2}\theta_k$, when the gain matrices $G_1$ and $G$ can be chosen freely, is the same for the two classes of stochastic iterations described above.*

PROOF. The single-time-scale iteration is of the form

$$Z_{k+1} = Z_k + \beta_k G(b - AZ_k + U_k),$$

where

$$Z_k = \begin{bmatrix} \theta_k \\ r_k \end{bmatrix}, \qquad U_k = \begin{bmatrix} V_k \\ W_k \end{bmatrix}$$

and

$$b = \begin{bmatrix} b_1 \\ b_2 \end{bmatrix}, \qquad A = \begin{bmatrix} A_{11} & A_{12} \\ A_{21} & A_{22} \end{bmatrix}.$$

As is well known [Kushner and Yin (1997)], the optimal (in the sense of positive definiteness) asymptotic covariance of $\beta_k^{-1/2} Z_k$ over all possible choices of $G$ is the covariance of $A^{-1} U_k$. We note that the top block of $A^{-1} U_k$ is equal to $\Delta^{-1}(V_k - A_{12} A_{22}^{-1} W_k)$. It then follows that the optimal asymptotic covariance matrix of $\beta_k^{-1/2} \theta_k$ is the covariance of $\Delta^{-1}(V_k - A_{12} A_{22}^{-1} W_k)$.

For the two-time-scale iteration, Theorem 2.8 shows that for any choice of $G_1$, the asymptotic covariance is the same as for the single-time-scale iteration:

$$\theta_{k+1} = \theta_k + \beta_k G_1 (b_1 - \Delta \theta_k + V_k - A_{12} A_{22}^{-1} W_k).$$

From this, it follows that the optimal asymptotic covariance of $\beta_k^{-1/2} \theta_k$ is the covariance of $\Delta^{-1}(V_k - A_{12} A_{22}^{-1} W_k)$, which is the same as for single-time-scale iterations. □

**4. Asymptotic normality.** In Section 2, we showed that $\beta_k^{-1} E[\hat{\theta}_k \hat{\theta}_k']$ converges to $\Sigma_{11}^{(0)}$. The proof techniques used in that section do not extend easily (without stronger assumptions) to the nonlinear case. For this reason, we develop here a different result, namely, the asymptotic normality of $\hat{\theta}_k$, which is easier to extend to the nonlinear case. In particular, we show that the distribution of $\beta_k^{-1/2} \hat{\theta}_k$ converges to a zero-mean normal distribution with covariance matrix $\Sigma_{11}^{(0)}$. The proof is similar to the one presented in Polyak (1990) for stochastic approximation with averaging.

THEOREM 4.1. *If Assumptions 2.1–2.5 hold with $\varepsilon = 0$, then $\beta_k^{-1/2} \hat{\theta}_k$ converges in distribution to $N(0, \Sigma_{11}^{(0)})$.*

PROOF. Recall the iterations (2.11) in terms of transformed variables $\tilde{\theta}$ and $\tilde{r}$. Assuming that $k$ is large enough so that $B_{21}^k = 0$, these iterations can be written as

$$\tilde{\theta}_{k+1} = (I - \beta_k \Delta) \tilde{\theta}_k - \beta_k A_{12} \tilde{r}_k + \beta_k V_k + \beta_k \delta_k^{(1)},$$

$$\tilde{r}_{k+1} = (I - \gamma_k A_{22}) \tilde{r}_k + \gamma_k W_k + \beta_k \delta_k^{(2)} + \beta_k (L_{k+1} + A_{22}^{-1} A_{21}) V_k,$$

where $\delta_k^{(1)}$ and $\delta_k^{(2)}$ are given by

$$\delta_k^{(1)} = A_{12} L_k \tilde{\theta}_k,$$

$$\delta_k^{(2)} = -(L_{k+1} + A_{22}^{-1} A_{21}) A_{12} \tilde{r}_k.$$

Using Theorem 2.6, $E[|\tilde{\theta}_k|^2]/\beta_k$ and $E[|\tilde{r}_k|^2]/\gamma_k$ are bounded, which implies that

$$E[|\delta_k^{(1)}|^2] \leq c \beta_k |L_k|^2,$$

(4.1)
$$E[|\delta_k^{(2)}|^2] \leq c\gamma_k,$$

for some constant $c > 0$. Without loss of generality assume $k_0 = 0$ in (2.11). For each $i$, define the sequence of matrices $\Theta_j^i$ and $R_j^i$, $j \geq i$, as

$$\begin{aligned}\Theta_i^i &= I, \\ \Theta_{j+1}^i &= \Theta_j^i - \beta_j \Delta \Theta_j^i \quad \forall j \geq i, \\ R_i^i &= I, \\ R_{j+1}^i &= R_j^i - \gamma_j A_{22} R_j^i \quad \forall j \geq i.\end{aligned}$$

Using the above matrices, $\tilde{r}_k$ and $\tilde{\theta}_k$ can be rewritten as

(4.2) $$\tilde{\theta}_k = \Theta_k^0 \tilde{\theta}_0 - \sum_{i=0}^{k-1} \beta_i \Theta_k^i A_{12} \tilde{r}_i + \sum_{i=0}^{k-1} \beta_i \Theta_k^i V_i + \sum_{i=0}^{k-1} \beta_i \Theta_k^i \delta_i^{(1)}$$

and

(4.3)
$$\begin{aligned}\tilde{r}_k = R_k^0 \tilde{r}_0 &+ \sum_{i=0}^{k-1} \gamma_i R_k^i W_i + \sum_{i=0}^{k-1} \beta_i R_k^i \delta_i^{(2)} \\ &+ \sum_{i=0}^{k-1} \beta_i R_k^i (L_{i+1} + A_{22}^{-1} A_{21}) V_i.\end{aligned}$$

Substituting the right-hand side of (4.3) for $\tilde{r}_k$ in (4.2), and dividing by $\beta_k^{1/2}$, we have

(4.4)
$$\begin{aligned}\beta_k^{-1/2} \tilde{\theta}_k = \frac{1}{\sqrt{\beta_0}} \tilde{\Theta}_k^0 \tilde{\theta}_0 &+ \sum_{i=0}^{k-1} \beta_i \tilde{\Theta}_k^i A_{12} (\beta_i^{-1/2} R_i^0 \tilde{r}_0) \\ &+ \sum_{i=0}^{k-1} \beta_i \tilde{\Theta}_k^i (\beta_i^{-1/2} \delta_i^{(1)}) + S_k^{(1)} + S_k^{(2)} + S_k^{(3)} \\ &+ \sum_{i=0}^{k-1} \sqrt{\beta_i} \tilde{\Theta}_k^i (V_i + A_{12} A_{22}^{-1} W_i),\end{aligned}$$

where

$$\tilde{\Theta}_k^i = \sqrt{\frac{\beta_i}{\beta_k}} \Theta_k^i \quad \forall k \geq i,$$

$$S_k^{(1)} = \sum_{i=0}^{k-1} \beta_i \tilde{\Theta}_k^i A_{12} \left( \beta_i^{-1/2} \sum_{j=0}^{i-1} \beta_j R_i^j \delta_j^{(2)} \right),$$

$$S_k^{(2)} = \sum_{i=0}^{k-1} \beta_i \tilde{\Theta}_k^i A_{12} \left( \beta_i^{-1/2} \sum_{j=0}^{i-1} \beta_j R_i^j (L_{j+1} + A_{22}^{-1} A_{21}) V_j \right),$$

$$S_k^{(3)} = \sum_{i=0}^{k-1} \sqrt{\beta_i} \tilde{\Theta}_k^i A_{12} \sum_{j=0}^{i-1} \gamma_j R_i^j W_j - \sum_{j=0}^{k-1} \sqrt{\beta_j} \tilde{\Theta}_k^j A_{12} A_{22}^{-1} W_j.$$

We wish to prove that the various terms in (4.4), with the exception of the last one, converge in probability to zero. Note that the last term is a martingale and therefore, can be handled by appealing to a central limit theorem for martingales. Some of the issues we encounter in the remainder of the proof are quite standard, and in such cases we will only provide an outline.

To better handle each of the various terms in (4.4), we need approximations of $\Theta_k^i$ and $R_k^i$. To do this, consider the nonlinear map $A \mapsto \exp(A)$ from square matrices to square matrices. A simple application of the inverse function theorem shows that this map is a diffeomorphism (differentiable, one-to-one with differentiable inverse) in a neighborhood of the origin. Let us denote the inverse of $\exp(\cdot)$ by $\ln(\cdot)$. Since $\ln(\cdot)$ is differentiable around $I = \exp(0)$, the function $\varepsilon \mapsto \ln(I - \varepsilon A)$ can be expanded into Taylor's series for sufficiently small $\varepsilon$ as follows:

$$\ln(I - \varepsilon A) = -\varepsilon(A - E(\varepsilon)),$$

where $E(\varepsilon)$ commutes with $A$ and $\lim_{\varepsilon \to 0} E(\varepsilon) = 0$. Assuming, without loss of generality, that $\gamma_0$ and $\beta_0$ are small enough for the above approximation to hold, we have for $k \geq 0$,

(4.5)
$$\Theta_k^i = \exp\left( -\sum_{j=i}^{k-1} \beta_j (\Delta - E_j^{(1)}) \right),$$

$$R_k^i = \exp\left( -\sum_{j=i}^{k-1} \gamma_j (A_{22} - E_j^{(2)}) \right),$$

for some sequence of matrices $\{E_k^{(i)}\}$, $i = 1, 2$, converging to zero. To obtain a similar representation for $\tilde{\Theta}_k^i$, note that Assumption 2.5(1) implies

(4.6) $$\frac{\beta_k}{\beta_{k+1}} = (1 + \beta_k(\varepsilon_k + \bar{\beta})),$$

for some $\varepsilon_k \to 0$. Therefore, using the fact that $1 + x = \exp(x(1 - o(x)))$ and (4.5), we have

(4.7) $$\tilde{\Theta}_k^i = \exp\left( -\sum_{j=i}^{k-1} \beta_j \left( \left( \Delta - \frac{\bar{\beta}}{2} I \right) - E_j^{(3)} \right) \right),$$

for some sequences of matrices $E_k^{(3)}$ converging to zero. Furthermore, it is not difficult to see that the matrices $E_k^{(i)}$, $i = 1, 2, 3$, commute with the matrices $\Delta$, $A_{22}$ and $\Delta - (\bar{\beta}/2)I$, respectively. Since $-\Delta$, $-(\Delta - (\bar{\beta}/2)I)$ and $-A_{22}$ are Hurwitz, using standard Lyapunov techniques we have for some constants $c_1, c_2 > 0$,

(4.8)
$$\max(|\Theta_k^i|, |\tilde{\Theta}_k^i|) \leq c_1 \exp\left(-c_2 \sum_{j=i}^{k-1} \beta_j\right),$$

$$|R_k^i| \leq c_1 \exp\left(-c_2 \sum_{j=i}^{k-1} \gamma_j\right).$$

Therefore it is easy to see that the first term in (4.4) goes to zero w.p.1. To prove that the second term goes to zero w.p.1, note that $\ln \beta_i \approx -\bar{\beta} \sum_{j=0}^{i-1} \beta_j$ [cf. (4.6)] and therefore for some $c_1, c_2 > 0$,

$$|\beta_i^{-1/2} R_i^0 \tilde{r}_0| \leq c_1 \exp\left(-c_2 \sum_{j=0}^{i-1}\left(\gamma_j - \frac{\bar{\beta}}{2}\beta_j\right)\right),$$

which goes to zero as $i \to \infty$ (Assumption 2.3). Therefore, it follows from Lemma A.3 that the second term also converges to zero w.p.1. Using (4.1) and Lemma A.3, it is easy to see that the third term in (4.4) converges in the mean (i.e., in $L_1$) to zero. Next, consider $E[|S_k^{(1)}|]$. Using (4.1), we have for some positive constants $c_1, c_2$ and $c_3$,

$$E\left[\left|\beta_i^{-1/2}\sum_{j=0}^{i-1}\beta_j R_j^i \delta_j^{(2)}\right|\right]$$
$$\leq c_1 \sum_{j=0}^{i-1} \gamma_j \exp\left(-\sum_{l=j}^{i-1}(c_2 \gamma_l - c_3 \beta_l)\right)\sqrt{\frac{\beta_j}{\gamma_j}}.$$

Since $\beta_j/\gamma_j \to 0$, Lemma A.3 implies that $S_k^{(1)}$ converges in the mean to zero. To study $S_k^{(2)}$, consider

$$E\left[\left|\beta_i^{-1/2}\sum_{j=0}^{i-1}\beta_j R_i^j (L_{j+1} + A_{22}^{-1}A_{21})V_j\right|^2\right].$$

Since the $V_k$ are zero mean i.i.d., the above term is bounded above by

$$c_1 \sum_{j=0}^{i-1} \gamma_j \exp\left(-\sum_{l=j}^{i-1}(c_2\gamma_l - c_3\beta_l)\right)\frac{\beta_j}{\gamma_j}$$

for some constants $c_1, c_2$ and $c_3$. Lemma A.3 implies that $S_k^{(2)}$ converges in the mean to zero. Finally, consider $S_k^{(3)}$. By interchanging the order of summation, it can be rewritten as

$$\text{(4.9)} \qquad \sum_{j=0}^{k-1} \sqrt{\beta_j} \tilde{\Theta}_k^j \left[ \frac{\gamma_j}{\beta_j} \sum_{i=j}^{k-1} \beta_i (\Theta_i^j)^{-1} A_{12} R_i^j - A_{12} A_{22}^{-1} \right] W_j.$$

Since $-A_{22}$ is Hurwitz, we have

$$A_{22}^{-1} = \int_0^\infty \exp(-A_{22} t) \, dt,$$

and we can rewrite the term inside the brackets in (4.9) as

$$\sum_{i=j}^{k-1} \gamma_i \left( \frac{\gamma_j \beta_i}{\beta_j \gamma_i} (\Theta_i^j)^{-1} - I \right) A_{12} R_i^j$$

$$+ A_{12} \left( \sum_{i=j}^{k-1} \gamma_i R_i^j - \int_0^{\sum_{i=j}^{k-1} \gamma_i} \exp(-A_{22} t) \, dt \right) - A_{12} A_{22}^{-1} \exp\left( -\sum_{i=j}^{k-1} \gamma_i A_{22} \right).$$

We consider each of these terms separately. To analyze the first term, we wish to obtain an "exponential" representation for $\gamma_j \beta_i / \beta_j \gamma_i$. It is not difficult to see from Assumptions 2.5 (1) and (2) that

$$\frac{\beta_{k+1}}{\gamma_{k+1}} = \frac{\beta_k}{\gamma_k}(1 - \varepsilon_k \gamma_k)$$

$$= \frac{\beta_k}{\gamma_k} \exp(-\varepsilon_k \gamma_k + O(\varepsilon_k^2 \gamma_k^2)),$$

where $\varepsilon_k \to 0$. Therefore, using (4.5) and the mean value theorem, we have

$$\left| \frac{\gamma_j \beta_i}{\beta_j \gamma_i} (\Theta_i^j)^{-1} - I \right|$$

$$\leq c_1 \sup_{l \geq j} \left( \varepsilon_l + \frac{\beta_l}{\gamma_l} \right) \left( \sum_{l=j}^{i-1} \gamma_l \right) \exp\left( c_2 \sum_{l=j}^{i-1} \left( \varepsilon_l + \frac{\beta_l}{\gamma_l} \right) \gamma_l \right),$$

which in turn implies, along with Lemma A.4 (with $p = 1$) and Assumption 2.3, that the first term is bounded in norm by $c \sup_{l \geq j}(\varepsilon_l + \gamma_l/\beta_l)$ for some constant $c > 0$. The second term is the difference between an integral and its Riemannian approximation and therefore is bounded in norm by $c \sup_{l \geq j} \gamma_l$ for some constant $c > 0$. Finally, since $-A_{22}$ is Hurwitz, the norm of the third term is bounded above by

$$c_1 \exp\left( -c_2 \sum_{i=j}^{k-1} \gamma_i \right)$$

for some constants $c_1, c_2 > 0$. An explicit computation of $E[|S_k^{(3)}|^2]$, using the fact that $(V_k, W_k)$ is zero-mean i.i.d., and an application of Lemma A.3 shows that $S_k^{(3)}$ converges to zero in the mean square. Therefore, the distribution of $\beta_k^{-1/2} \tilde{\theta}_k$ converges to the asymptotic distribution of the martingale comprising the remaining terms. To complete the proof, we use the standard central limit theorem for martingales [see Duflo (1997)]. The key assumption of this theorem is Lindberg's condition which, in our case, boils down to the following: for each $\varepsilon > 0$,

$$\lim_k \sum_{i=0}^{k-1} E\Big[|X_i^{(k)}|^2 I\{|X_i^{(k)}| \geq \varepsilon\}\Big] = 0,$$

where $I$ is the indicator function and for each $i < k$,

$$X_i^{(k)} = \sqrt{\beta_i} \tilde{\Theta}_k^i (V_i + A_{12} A_{22}^{-1} W_i).$$

The verification of this assumption is quite standard. $\square$

REMARK. Similar results are possible for nonlinear iterations with Markov noise. For an informal sketch of such results, see Konda (2002).

## APPENDIX: AUXILIARY RESULTS

**A.1. Verification of (2.11).** Without loss of generality, assume that $b_1 = b_2 = 0$. Then, $\theta^* = 0$ and

$$\tilde{\theta}_k = \hat{\theta}_k = \theta_k,$$

and, using the definition of $\tilde{r}_k$ [cf. (2.4) and (2.10)], we have

(A.10) $\quad \tilde{r}_k = L_k \theta_k + \hat{r}_k = L_k \theta_k + r_k + A_{22}^{-1} A_{21} \theta_k = r_k + M_k \theta_k,$

where

$$M_k = L_k + A_{22}^{-1} A_{21}.$$

To verify the equation for $\tilde{\theta}_{k+1} = \theta_{k+1}$, we use the recursion for $\theta_{k+1}$, to obtain

$$\begin{aligned}
\theta_{k+1} &= \theta_k - \beta_k(A_{11}\theta_k + A_{12}r_k - V_k) \\
&= \theta_k - \beta_k(A_{11}\theta_k + A_{12}\tilde{r}_k - A_{12}(L_k + A_{22}^{-1}A_{21})\theta_k - V_k) \\
&= \theta_k - \beta_k(A_{11}\theta_k - A_{12}A_{22}^{-1}A_{21}\theta_k - A_{12}L_k\theta_k + A_{12}\tilde{r}_k - V_k) \\
&= \theta_k - \beta_k(\Delta\theta_k - A_{12}L_k\theta_k + A_{12}\tilde{r}_k) + \beta_k V_k \\
&= \theta_k - \beta_k(B_{11}^k \theta_k + A_{12}\tilde{r}_k) + \beta_k V_k,
\end{aligned}$$

where the last step makes use of the definition $B_{11}^k = \Delta - A_{12}L_k$.

To verify the equation for $\tilde{r}_{k+1}$, we first use the definition (A.10) of $\tilde{r}_{k+1}$, and then the update formulas for $\theta_{k+1}$ and $r_{k+1}$, to obtain

$$
\begin{aligned}
\tilde{r}_{k+1} &= r_{k+1} + (A_{22}^{-1}A_{21} + L_{k+1})\theta_{k+1} \\
&= r_k - \gamma_k(A_{21}\theta_k + A_{22}r_k - W_k) + (A_{22}^{-1}A_{21} + L_{k+1})\theta_{k+1} \\
&= r_k - \gamma_k(A_{21}\theta_k + A_{22}(\tilde{r}_k - (L_k + A_{22}^{-1}A_{21})\theta_k) - W_k) \\
&\quad + (A_{22}^{-1}A_{21} + L_{k+1})\theta_{k+1} \\
&= r_k - \gamma_k(A_{22}\tilde{r}_k - A_{22}L_k\theta_k - W_k) + M_{k+1}\theta_{k+1} \\
&= r_k + M_{k+1}\theta_k - \gamma_k(A_{22}\tilde{r}_k - A_{22}L_k\theta_k - W_k) \\
&\quad - \beta_k M_{k+1}(B_{11}^k\theta_k + A_{12}\tilde{r}_k - V_k) \\
&= r_k + M_k\theta_k - \gamma_k\left[\frac{L_k - L_{k+1}}{\gamma_k} - A_{22}L_k + \frac{\beta_k}{\gamma_k}M_{k+1}B_{11}^k\right]\theta_k \\
&\quad + \gamma_k W_k - \gamma_k\left(A_{22} + \frac{\beta_k}{\gamma_k}M_{k+1}A_{12}\right)\tilde{r}_k + \beta_k M_{k+1}V_k \\
&= \tilde{r}_k - \gamma_k(B_{21}^k\tilde{\theta}_k + B_{22}^k\tilde{r}_k) + \gamma_k W_k + \beta_k M_{k+1}V_k,
\end{aligned}
$$

which is the desired formula.

### A.2. Convergence of the recursion (2.12).

LEMMA A.1.  *For $k_0$ sufficiently large, the (deterministic) sequence of matrices $\{L_k\}$ defined by (2.12) is well defined and converges to zero.*

PROOF.  The recursion (2.12) can be rewritten, for $k \geq k_0$, as

(A.2)
$$
\begin{aligned}
L_{k+1} &= (I - \gamma_k A_{22})L_k \\
&\quad + \beta_k(A_{22}^{-1}A_{21}B_{11}^k + (I - \gamma_k A_{22})L_k B_{11}^k)(I - \beta_k B_{11}^k)^{-1},
\end{aligned}
$$

which is of the form

$$L_{k+1} = (I - \gamma_k A_{22})L_k + \beta_k D_k(L_k),$$

for a sequence of matrix-valued functions $D_k(L_k)$ defined in the obvious manner. Since $-A_{22}$ is Hurwitz, there exists a quadratic norm

$$|x|_Q = \sqrt{x'Qx},$$

a corresponding induced matrix norm, and a constant $a > 0$ such that

$$|(I - \gamma A_{22})|_Q \leq (1 - a\gamma)$$

for every sufficiently small $\gamma$. It follows that

$$|(I - \gamma A_{22})L|_Q \leq (1 - a\gamma)|L|_Q$$

for all matrices $L$ of appropriate dimensions and for $\gamma$ sufficiently small. Therefore, for sufficiently large $k$, we have

$$|L_{k+1}|_Q \leq (1 - \gamma_k a)|L_k|_Q + \beta_k |D(L_k)|_Q.$$

For $k_0$ sufficiently large, the sequence of functions $\{D_k(\cdot)\}_{k \geq k_0}$ is well defined and uniformly bounded on the unit $Q$-ball $\{L : |L|_Q \leq 1\}$. To see this, note that as long as $|L_k|_Q \leq 1$, we have $|B_{11}^k| = |\Delta - A_{12}L_k| \leq c$, for some absolute constant $c$. With $\beta_k$ small enough, the matrix $I - \beta_k B_{11}^k$ is invertible, and satisfies $|(I - \beta_k B_{11}^k)^{-1}| \leq 2$. With $|B_{11}^k|$ bounded by $c$, we have

$$|A_{22}^{-1} A_{21} B_{11}^k + (I - \gamma_k A_{22}) L_k B_{11}^k| \leq d(1 + |L_k|),$$

for some absolute constant $d$. To summarize, for large $k$, if $|L_k|_Q \leq 1$, we have $|D_k(L_k)| \leq 4d$. Since any two norms on a finite-dimensional vector space are equivalent, we have

$$|L_{k+1}|_Q \leq (1 - \gamma_k a)|L_k|_Q + (\gamma_k a)\left(\frac{d_1 \beta_k}{a \gamma_k}\right),$$

for some constant $d_1 > 0$. Recall now that the sequence $L_k$ is initialized with $L_{k_0} = 0$. If $k_0$ is large enough so that $d_1 \beta_k / a \gamma_k < 1$, then $|L_k|_Q \leq 1$ for all $k$. Furthermore, since $1 - x \leq e^{-x}$, we have

$$|L_k|_Q \leq \sum_{j=k_0}^{k-1} \gamma_j \exp\left(-a \sum_{i=j}^{k-1} \gamma_i\right)\left(\frac{d_1 \beta_j}{\gamma_j}\right).$$

The rest follows from Lemma A.3 as $\beta_k / \gamma_k \to 0$. □

**A.3. Linear matrix iterations.** Consider a linear matrix iteration of the form

$$\Sigma_{k+1} = \Sigma_k + \beta_k(\Gamma - A\Sigma_k - \Sigma_k B + \delta_k(\Sigma_k))$$

for some square matrices $A$, $B$, step-size sequence $\beta_k$ and sequence of matrix-valued affine functions $\delta_k(\cdot)$. Assume:

1. The real parts of the eigenvalues of $A$ are positive and the real parts of the eigenvalues of $B$ are nonnegative. (The roles of $A$ and $B$ can also be interchanged.)
2. $\beta_k$ is positive and

$$\beta_k \to 0, \qquad \sum_k \beta_k = \infty.$$

3. $\lim_k \delta_k(\cdot) = 0$.

We then have the following standard result whose proof can be found, for example, in Polyak (1976).

LEMMA A.2. *For any $\Sigma_0$, $\lim_k \Sigma_k = \Sigma^*$ exists and is the unique solution to the equation*

$$A\Sigma + \Sigma B = \Gamma.$$

**A.4. Convergence of some series.** We provide here some lemmas that are used in the proof of asymptotic normality. Throughout this section, $\{\gamma_k\}$ is a positive sequence such that:

1. $\gamma_k \to 0$, and
2. $\sum_k \gamma_k = \infty$.

Furthermore, $\{t_k\}$ is the sequence defined by

$$t_0 = 0, \qquad t_k = \sum_{j=0}^{k-1} \gamma_k, \qquad k > 0.$$

LEMMA A.3. *For any nonnegative sequence $\{\delta_k\}$ that converges to zero and any $p \geq 0$, we have*

(A.3) $$\lim_k \sum_{j=0}^{k} \gamma_j \left( \sum_{i=j}^{k-1} \gamma_i \right)^p \exp\left( -\sum_{i=j}^{k-1} \gamma_i \right) \delta_j = 0.$$

PROOF. Let $\delta(\cdot)$ be a nonnegative function on $[0, \infty)$ defined by

$$\delta(t) = \delta_k, \qquad t_k \leq t < t_{k+1}.$$

Then it is easy to see that for any $k_0 > 0$,

$$\sum_{j=k_0}^{k} \gamma_j \left( \sum_{i=j}^{k-1} \gamma_i \right)^p \exp\left( -\sum_{i=j}^{k-1} \gamma_i \right) \delta_j$$
$$= \int_{t_{k_0}}^{t_k} (t_k - s)^p e^{-(t_k - s)} \delta(s)\, ds + e_k^{k_0},$$

where

$$|e_k^{k_0}| \leq c \sum_{j=k_0}^{k} \gamma_j^2 \left( \sum_{i=j}^{k-1} \gamma_i \right)^p \exp\left( -\sum_{i=j}^{k-1} \gamma_i \right) \delta_j$$

for some constant $c > 0$. Therefore, for $k_0$ sufficiently large, we have

$$\lim_k \sum_{j=k_0}^{k} \gamma_j \left( \sum_{i=j}^{k-1} \gamma_i \right)^p \exp\left( -\sum_{i=j}^{k-1} \gamma_i \right) \delta_j$$
$$\leq \frac{\lim_t \int_0^t \delta(s)(t-s)^p e^{-(t-s)}\, ds}{1 - c\sup_{k \geq k_0} \gamma_k}.$$

To calculate the above limit, note that

$$\lim_t \left| \int_0^t (t-s)^p e^{-(t-s)} \delta(s)\, ds \right|$$

$$= \lim_t \left| \int_0^t s^p e^{-s} \delta(t-s)\, ds \right|$$

$$\leq \lim_t \left( \sup_{s \geq t-T} |\delta(s)| \right) \int_0^T s^p e^{-s}\, ds + \sup_s |\delta(s)| \int_T^\infty s^p e^{-s}\, ds$$

$$= \sup_s |\delta(s)| \int_T^\infty s^p e^{-s}\, ds.$$

Since $T$ is arbitrary, the above limit is zero. Finally, note that the limit in (A.3) does not depend on the starting limit of the summation. $\square$

LEMMA A.4. *For each $p \geq 0$, there exists $K_p > 0$ such that for any $k \geq j \geq 0$,*

$$\sum_{i=j}^k \gamma_i \left( \sum_{l=j}^{i-1} \gamma_l \right)^p \exp\left( -\sum_{l=j}^{i-1} \gamma_l \right) \leq K_p.$$

PROOF. For all $j$ sufficiently large, we have

$$\sum_{i=j}^k \gamma_i \left( \sum_{l=j}^{i-1} \gamma_l \right)^p \exp\left( -\sum_{l=j}^{i-1} \gamma_l \right) \leq \frac{\int_0^{(t_k - t_j)} \tau^p e^{-\tau}\, d\tau}{1 - c \sup_{l \geq j} \gamma_l},$$

for some $c \geq 0$. $\square$

Laboratory for Information
and Decision Systems
Massachusetts Institute of Technology
77 Massachusetts Avenue
Cambridge, Massachusetts 02139
USA
e-mail: konda@alum.mit.edu
e-mail: jnt@mit.edu